\documentclass{amsart}
\usepackage{amscd,amssymb}

\theoremstyle{plain}
\newtheorem{thm}{Theorem}
\newtheorem{lem}[thm]{Lemma}
\newtheorem{definicja}[thm]{Definition}
\newtheorem{cor}[thm]{Corollary}
\newtheorem{prop}[thm]{Proposition}

\theoremstyle{definition}

\theoremstyle{remark}

\numberwithin{equation}{section}

\usepackage{pst-all}

\begin{document}
\setcounter{section}{-1}

\title{Local operations and eventually open actions 
} 

\author{Barbara Majcher-Iwanow} 

\maketitle 
\centerline{Institute of Mathematics, Wroc{\l}aw University, pl.Grunwaldzki 2/4, }
\centerline{50-384, Wroc{\l}aw, Poland}  
\centerline{{\em E-mail address}:  biwanow@math.uni.wroc.pl ;  
{\em Fax number}: 48-71-3757429 }

%
 
\begin{abstract} 
We study continuous actions of Polish groups on Polish spaces.
We develop Scott analysis introduced by Hjorth for studying orbit equivalence relations.
We define {\em eventually open} actions and prove that this 
property characterizes the actions endowed with a complete system of
hereditarily countable invariant structures. \parskip0pt 

{\em Keywords}: Polish G-spaces; Scott analysis; Canonical partitions. 
\parskip0pt 

{\em 2000 Mathematics Subject Classification}: 03E15.  
\end{abstract} 

\topskip 20pt
\section{Introduction}

In this paper we study continuous actions of Polish groups
on Polish spaces (say $G$ on $X$) by means of generalized Scott
invariants introduced by Hjorth in \cite{hjorth}.
Modifying the generalized Hjorth-Scott analysis we
approach the orbit equivalence relation in a fashion
which exploits descriptive set theoretical view-point
slightly more intensively.\parskip0pt

The basic tool of Hjorth's work are
hereditarily countable structures $\phi_{\alpha}(x,U,V )$,
$U \subseteq_{open} X$, $V\subseteq_{open} G$,
corresponding to Scott characteristics (\cite{hjorth},  Chapter 6.2).
We rather concentrate on associated sets
$\  B_{\alpha}(x,U,V)=\{y\in X: \phi_{\alpha}(y,U,V )=\phi_{\alpha}(x,U,V)\}$
(we call them $\alpha$-pieces) and their presentations with
use of some {\em operation of local saturation}.
This direction can be considered as a generalization
of the notion of {\em canonical partitions} introduced by Becker in \cite{becker}.
Following this way we are able to supplement Hjorth's work with
a couple of new statements concerning the sets $B_{\alpha}(x,U,V)$.
In particular in  Theorem \ref{VB} we present
a canonical form  for $B_{\alpha}(x,U,V)$.
This immediately implies that
the sets $B_{\alpha}(x,U,V)$ are Borel and
moreover this describes their Borel complexity.
\parskip0pt

The original motivation for this result is connected
with the problems of coding of $G$-orbits in admissible sets.
In order to extend the results of \cite{nadel} and \cite{basadm}
to the general case of Polish $G$-spaces, Theorem \ref{VB}
looks very helpful.
This stuff will be considered in a separate paper.

In this paper we first concentrate on
refinnig topologies by extending the initial
basis by families of sets of the form
$B_{\beta}(x,U,V)$.
Applying Theorem \ref{VB} we show that
the original topology enriched upon
some natural families of these sets
generates on $X$ finer topologies endowed with
the same Borel structure as the initial one so that
each $B_{\alpha}(x,X,G)$ with the corresponding subspace
topology becomes a Polish $G$-space.
This generalizes a similar theorem proved by
Hjorth in \cite{hjorth} in the case when $\alpha$ is
$\gamma^{\star}(x)$, the generalized Scott rank of $x$.

We then study when the maps $G\rightarrow Gx$
defined by $g\rightarrow gx$ are open (for all $x\in X$)
with respect to appropriate topologies mentioned above.
We prove that this property is equivalent to
the property that the generalized Scott analysis
applied to the orbit equivalence relation
leads to a complete system of invariants.
Moreover we can restate this
as a very simple condition which we
call {\em eventual openness of the action}
(this is the content of Theorem \ref{claschar}).

It is worth noting that it is proved in \cite{hjorth} that
orbit equivalence relations equipped with a complete system of
generalized Scott invariants are classifiable
by countable models, i.e. they are Borel reducible to
the isomorphism relation on the space $Mod(\mathcal{L})$
of all countable structures of some countable language $\mathcal{L}$.
\parskip0pt

The operation of {\em local saturation} and its basic
properties are presented in Section 1.
Section 2 is devoted to our approach to the generalized
Hjorth-Scott analysis.
Eventually open actions are studied in Section 3.
Along with the local counterpart of saturation we apply there
a local version of Vaught transforms.
This may be interesting in itself.

\section{Preliminaries}

In the first part of this section we recall standard notation and
facts concerning Polish group actions.
In the second one we define {\em local saturation} -
a new operation arising in this context.
This operation  is of particular importance for this paper.

\subsection{Notation.}

A {\em Polish space} ({\em group}) is a separable, completely
metrizable topological space (group).
If a Polish group $G$ continuously acts on a Polish space $X$,
then we say that $X$ is a {\em Polish $G$-space}.
We  say that a subset of $X$ is {\em invariant} if
it is $G$-invariant.
All basic facts concerning Polish $G$-spaces can
be found in \cite{bk}, \cite{hjorth} and \cite{kechris}.\parskip0pt

Let $G$ be a Polish group,
${\mathcal N}$ be a countable basis of $G$ and $G_0=\{g_i:i\in \omega\}$
be a countable dense subgroup of $G$.
Let ${\mathcal V}\subseteq {\mathcal N}$,
${\mathcal V}=\{V_m: m\in \omega\}$, be a countable basis of open
neighborhoods of the unity of $G$.
We shall assume that $V=V^{-1}$ and $V^g\in {\mathcal V}$,
whenever $V\in {\mathcal V}$ and  $g\in G_0$.
Besides ${\mathcal V}$ we shall use the symbol $\hat{\mathcal V}$
to denote the set of all (not only basic)  symmetric neighourhoods
of the unity $1_G$.
\parskip0pt

Let $\langle X, \tau\rangle$ be a Polish $G$-space and
${\mathcal U}=\{U_n:n\in\omega\}$ be a countable basis of $X$.
We assume  that for every $U\in {\mathcal U}$ and $g\in G_0$ we have
$gU\in {\mathcal U}$.

Since we shall use Vaught transforms,
recall the corresponding definitions.
The Vaught $*$-transform of a set $B\subseteq X$
with respect to an open $H\subseteq G$ is the set
$B^{*H}=\{ x\in X:\{ g\in H:gx\in B\}$ is comeager in $H\}$,
the Vaught $\Delta$-transform of $B$ is the set
$B^{\Delta H}=\{ x\in X:\{ g\in H:gx\in B\}$
is not meager in $H\}$.
It is known that for any $x\in X$ and $g\in G$,
$gx\in B^{*H} \Leftrightarrow x\in B^{*Hg}$ and
$gx\in B^{\Delta H} \Leftrightarrow x\in B^{\Delta Hg}$.\parskip0pt

It is worth noting that
{\em for any open $B\subseteq X$ and any open $K<G$ we have}
$B^{\Delta K}=KB$, where $KB=\{gx:g\in K, x\in B\}$.
Indeed, by continuity of the action for any $x\in KB$ and
$g\in K$ with $gx\in B$ there are open neighbourhoods
$K_1\subseteq K$ and $B_1\subseteq KB$ of $g$ and $x$
respectively so that $K_1 B_1\subseteq B$;
thus $x\in B^{\Delta K}$.
Moreover for every countable ordinal $\alpha$ 
if $B\in \Sigma^{0}_{\alpha}(X)$,
then $B^{\Delta H}\in \Sigma^{0}_{\alpha}(X)$
and if $B\in\Pi^{0}_{\alpha}(X)$, then
$B^{*H}\in\Pi^{0}_{\alpha}(X)$.\parskip0pt

Other basic properties of Vaught transforms can be found
in \cite{bk} and \cite{kechris}. \parskip0pt

\subsection{Local saturation}

When $G$ admits a basis of open subgroups at its unity $1_G$,
then every Polish $G$-space admits a basis consisting of the sets
which are invariant with respect to some open basic subgroup
of $G$.
Since such a $G$ is isomorphic to a closed permutation
group (see \cite{bk} for details), we may easily generalize Scott method
described in 6.1 of \cite{hjorth} to analyse orbit equivalence
relations arising in these situations.\parskip0pt

To handle with difficulties of the general case we introduce
a local variant of the operation of saturation
which generalizes the concept of a {\em local orbit}
introduced in \cite{hjorth}.

\begin{definicja}\label{sat}
Let $U\subseteq X$ be open\footnote{Actually we do not need to demand
that $U$ is open,
the definition makes sense for any $U$.} and $V\in \hat{\mathcal V}$.
For every $A\subseteq X$
we define inductively an increasing sequence $(V^{[n]}_UA)_{n\in \omega}$
of subsets of $X$ as follows:
$$
\left\{\begin{array}{l@{\ = \ }l}
V^{[0]}_UA&A\cap U, \\ 
V^{[n+1]}_UA&V(V^{[n]}_UA)\cap U.
\end{array}\right.
$$
The set $V_U A=\bigcup\limits_{n\in\omega}V^{[n]}_UA$
is called the local $V_U$-saturation of $A$.
For every $x\in X$ we shall write  $V_Ux$ instead of $\ V_U\{x\}$
and call this set the local $V_U$-orbit of $x$.
\end{definicja}

We see that $V_UA$ is contained in $U$ and contains $A\cap U$.
It is nonempty if and only if $A\cap U\not=\emptyset$,
in particular $V_Ux\not=\emptyset$ if and only if $x\in U$.
$V_UA$ is a union of all local $V_U$-orbits of  elements of $A$.\parskip0pt

There is a natural correspondence
between local orbits and suitably defined subsets of $G$
(which already appeared in \cite{hjorth}).
Let us introduce the corresponding definition and formulate basic facts.

\begin{definicja}\label{H} Let $U\subseteq X$ be open, $A\subseteq U$
and $V\in \hat{\mathcal V}$.
For every $x\in U$ we define an increasing sequence of subsets
of $G$ as follows.
$$
\begin{array}{l@{\ = \ }l}
{\langle V\rangle}^x_U(0)&
\{1_G\} ,\\
{\langle V\rangle}^x_U(n+1)&
\{gh:ghx\in U, h\in {\langle V\rangle}^x_U(n), g\in V\}.
\end{array}
$$
Then we put ${\langle V\rangle}_U^x
=\bigcup\limits_{n\in\omega}{\langle V\rangle}^x_U(n)$.\\ 
If $x\not\in U$, we put ${\langle V\rangle}_U^x=\emptyset$.\\
Finally we define
${\langle V\rangle}^A_U=\bigcap\{{\langle V\rangle}_U^x:x\in A\}$.
\end{definicja}

By this definition we see that $V_Ux={\langle V\rangle}_U^xx$.
The following simple lemma collects the very basic properties
of introduced sets.

\begin{lem}\label{bH} Let $U\subseteq X$ be open,
$V\in \hat{\mathcal V}$ and $x\in X$.
Then\parskip0pt

$(1)$ ${\langle V\rangle}^x_U(n+1)=\{gh:h\in {\langle V\rangle}^x_U(1),
g\in {\langle V\rangle}_U^{hx}(n)\}$,
for every $n\in\omega$.\parskip3pt

$(2)$ If $n,k \ge 0$ and $h\in {\langle V\rangle}^x_U(k)$, then
${\langle V\rangle}^{hx}_U(n)h\subseteq {\langle V\rangle}^x_U(n+k)$ and
${\langle V\rangle}^{x}_U(n)h^{-1}\subseteq {\langle V\rangle}^{hx}_U(n+k)$.
Consequently, ${\langle V\rangle}^{hx}_Uh={\langle V\rangle}^x_U$.
\parskip3pt

$(3)$ For every $h\in {\langle V\rangle}_U^x(1)$, $n\in \omega$,
$f\in {\langle V\rangle}_U^{hx}(n)$
there are open $U_x\subseteq U$, $W_h\subseteq {\langle V\rangle}_U^x(1)$
and $W_f\subseteq {\langle V\rangle}_U^{hx}(n)$
such that $x\in U_x$, $h\in W_h$, $f\in W_f$  and
$f'h'\in {\langle V\rangle}_U^y(n+1)$,
for every $y\in U_x$, $h'\in W_h$ and $f'\in W_f$.\parskip0pt

Consequently ${\langle V\rangle}^x_U(n), n\in\omega$,
and ${\langle V\rangle}_U^x$ are open.
\end{lem}

{\em Proof.} (1), (2) are easy consequences of the definition, then
(3) follows by induction from continuity of the action. $\Box$
\bigskip

{\bf Remark.} If $x\in U$, then the equality
${\langle V \rangle}^{hx}_Uh={\langle V \rangle}^x_U$
is not  true unless $h\in {\langle V\rangle}^x_U$.
Applying point (2) of the lemma above
we can easily check that for any $h,h'\in G$,
the sets ${\langle V \rangle}^{hx}_Uh$,
${\langle V \rangle}^{h'x}_Uh'$ are either equal
or disjoint.
Thus the family
$\{{\langle V\rangle}^{hx}_Uh: h\in G, hx\in U\}$
is a partition of the set
$\{h\in G: hx\in U\}$ into open sets.\footnote{In this form,
i.e. as classes of the appropriate equivalence relation, the sets we are
discussing appear in \cite{hjorth}.}
\bigskip

By classical results orbits of elements under continuous
(Borel) actions are Borel sets.
If we slightly modify the proof,
we see that this remains true for local orbits.

\begin{cor} Every local orbit under a continuous action  is a Borel set.
\end{cor}

{\em Proof}. Let $U\subseteq X$ be open, $V\in \hat{\mathcal V}$
and $x\in X$.
Put  $W=\{g\in G: gx\in V_Ux\}$.
We have $Wx=V_Ux$ and $W={\langle V\rangle}^x_UG_x$, where $G_x$
is the stabilizer of $x$.
By Lemma \ref{bH}(3) we see that $W$ is an open subset of $G$.
Let $T_x$ be a Borel transversal of $G/G_x$.
Then $W\cap T_x$ is also Borel and the function
$g\to gx$ is a bijection from $W\cap T_x$ onto $V_Ux$.
Hence the latter has to be Borel.
$\Box$
\bigskip

The other simple properties of the operation of local saturation
are collected below. \parskip0pt

\begin{lem}\label{locsat}
Let  $U,U'\subseteq X$ be open and $V,V'\in \hat{\mathcal V}$.
For every $A,B,A_1,A_2,\ldots \subseteq X$, $x,y\in U$ and $f\in G$
the following statements hold.

$(1)$ If $A\subseteq B$, $U'\subseteq U$ and $V'\subseteq V$ then
${V'}_{U'}A\subseteq V_UB$.\parskip2pt

$(2)$ $V_U(\bigcup\limits_n A_n)=\bigcup\limits_n V_U A_n$.\parskip2pt

$(3)$ $V_U(V_UA)=V_UA$.\parskip2pt

$(4)$ $V_Ux=V_Uy$ if and only if $x\in V_Uy$.\parskip2pt

$(5)$ If $x\in U'$ and $V'\subseteq V$, then
$V_Ux=U\cap \bigcup\{(gV'\cap V'g)x\cap gU':\
g\in G_0\cap {\langle V\rangle}_U^x\}$.
\parskip2pt

$(6)$ $f(V_UA)=(V^f)_{fU}(fA)$.
\end{lem}

{\em Proof.} (1) and (2) are immediate. 

(3) By (1), (2) we have
$V_U(V_UA)=V_U(\bigcup\limits_nV^{[n]}_UA)=
\bigcup\limits_nV_U(V^{[n]}_UA)$.
Let $n\in\omega$ be arbitrary.
It follows from the definition that
$V_U(V^{[n]}_UA)=\bigcup\limits_{i\ge n}V^{[i]}_UA$.
Since the family $\{V^{[n]}_UA\}_{n\in\omega}$\ is increasing, we see that
$V_UA=\bigcup\limits_{i\ge n}V^{[i]}_UA$ which completes the proof.
\parskip2pt

(4) $(\Rightarrow)$ is obvious.\parskip0pt

$(\Leftarrow)$ immediately follows by Lemma \ref{bH}(2).\parskip2pt

(5) To prove  $\subseteq$ consider an arbitrary $h\in {\langle V\rangle}^x_U$.
By Lemma \ref{bH}(3) we can find an open set
$W\subseteq {\langle V\rangle}^x_U$
such that
$h\in W$ and  $Wx\subseteq U$.
We may additionaly demand that $W^{-1}W\subseteq V'$
(since $V'$ is symmetric, then also $WW^{-1}\subseteq V'$),
$W^{-1}Wx\subseteq U'$ and $WW^{-1}x\subseteq U'$.
Then for an arbitrary $g\in G_0\cap W$ we have
$h\in gV'\cap V'g$ and $hx\in gU'$.

For the converse inclusion observe that
$V'gx\cap U\subseteq V_Ux$ whenever $gx\in V_Ux$.\parskip0pt

(6) follows from the fact that for every $n\in \omega$
we have $f(V^{[n]}_UA)=(V^f)^{[n]}_{fU}(fA)$.
The latter can be obtained by an easy inductive argument.
$\Box$\bigskip

The new concept of saturation entails a new concept of invariantness
- {\em local invariantness}.

\begin{definicja}
Let $U\subseteq X$ be open, $V\in \hat{\mathcal V}$ and $A\subseteq X$.
We say that  $A$ is  locally
$V_U$-invariant if \  $V_UA=A\cap U$.
\end{definicja}

{\bf Remark.} It follows that $A$ is locally $V_U$-invariant
if and only if $V_Ux\subseteq A$, for every $x\in A$.
Observe also that $A$ is locally $V_U$-invariant whenever
$V(A\cap U)\cap U =A\cap U$.
Indeed, the equality $V(A\cap U)\cap U=A\cap U$ implies by induction that
for every $n\in \omega$ we have $V^{[n]}_UA=A\cap U$ and thus
$V_UA=A\cap U$.
On the other hand we have
$A\cap U\ \subseteq V(A\cap U)\cap U\ \subseteq V_UA$.
Thus the equality $V_UA=A\cap U$ implies
that $V(A\cap U)\cap U=A\cap U$.\parskip0pt

Obviously every $A$ such that $A\supseteq U$ or
$U\cap A=\emptyset$ \ is locally $V_U$-invariant.
Moreover it can be justified by easy straightforward arguments that
the family of all locally $V_U$-invariant subsets of $X$ forms a complete
Boolean algebra.
By Lemma \ref{locsat} we see that  for every $A$,
$V_UA$ is a $V_U$-invariant set containing $A\cap U$.
If $A$ is open then $V_UA$ is open by Definition \ref{sat}.

\section{Sets arising in Polish group actions}

This section can be considered both as systematization and some
improvement of the material contained in Section 6.2 of \cite{hjorth}.
It is divided into two subsections.
In the first one we modify the generalized
Scott analysis developed by Hjorth.
The basic tool of Hjorth's work are
hereditarily countable structures $\phi_{\alpha}(x,U_n,V_n )$
corresponding to Scott characteristics.
We suggest slightly different approach and
concentrate on associated sets
$\  B_{\alpha}(x,U,V)=\{y\in X: \phi_{\alpha}(y,U,V )=\phi_{\alpha}(x,U,V)\}$.
We characterize the sets $B_{\alpha}(x,U,V)$ with use
of the operation of local saturation and study them
slightly further in order to present this material
in a complete form.

This direction can be considered as a generalization
of the notion of {\em canonical partitions} (see \cite{becker}).
Following this way we are able to supplement Hjorth's work with
a couple of new statements concerning the sets $B_{\alpha}(x,U,V)$.
The main result of this part is Theorem \ref{VB}
which makes possible to express the sets
$B_{\alpha}(x,U,V)$ in a canonical form.
This possiblity is of fundamental importance for our study.
In particular it enables us to prove that the sets
$B_{\alpha}(x,U,V)$ are Borel and describe their Borel complexity.
\parskip0pt

The second subsection is devoted to refinnig topologies
by extending the initial basis by families of
$\beta$-pieces, $\beta<\alpha$.
Applying Theorem \ref{VB} we show that
the original basis enriched upon the family
$\{B_{\beta}(x',U_n,V_m):x'\in V_Ux\cap U_n,\ n,m\in \omega,\ \beta<\alpha\}$
generates on $X$ a finer topology endowed with
the same Borel structure as the initial one so that
$B_{\alpha}(x,X,G)$ with the corresponding subspace
topology becomes a Polish $G$-space.
It is worth noting that a similar theorem is proved by
Hjorth in \cite{hjorth} in the case when $\alpha$ is
$\gamma^{\star}(x)$, the generalized Scott rank of $x$.
Thus our theorem can be considered as a generalization of it.

\subsection{Borel partitions}

As we have already mentioned $\alpha$-invariants
$\phi_{\alpha}(x,U,V)$ were introduced by Hjorth
as a counterpart of $\alpha$-invariants studied by Scott.
From now on we fix a countable basis
${\mathcal U}=\{U_n:n\in \omega\}$ of $X$
and a countable basis ${\mathcal V}=\{V_n:n\in \omega\}$
of open symmetric neighbourhoods of $1_G$.

\begin{definicja}\label{hj}(Hjorth)
For every  $U\in {\mathcal U}$, $x\in U$ and $V\in {\mathcal V}$
we define a set $\phi_{\alpha}(x,U,V)$ by
simultaneous induction on the ordinal $\alpha$: 
$$
\begin{array}{l@{\ \ }l}
\phi_1(x,U,V) & =\{ l:\ U_l\in {\mathcal U},
V_U x\cap U_l\not=\emptyset\}, \\ 
\phi_{\alpha +1}(x,U,V)&=
\{\langle\phi_{\alpha}(x',U_n,V_m), n, m\rangle:
x'\in V_Ux, U_n\subseteq U, V_m\subseteq  V\}, \\
\phi_{\lambda}(x,U,V) &=
\{\langle\phi_{\alpha}(x,U,V),\alpha\rangle:\alpha<\lambda\}
\mbox{ for }\lambda \mbox{ limit } .
\end{array}
$$
\end{definicja}

Every $\phi_{\alpha}(x,U,V)$, $1\le \alpha<\omega_1$,
defines the set $B_{\alpha}(x,U,V)=\{y\in U:
\phi_{\alpha}(y,U,V )=\phi_{\alpha}(x,U,V)\}$.
We call the sets of this form {\em $\alpha$-pieces}.
Additionally we treat every basic open $U$ as a 0-piece.
In the lemma below we put together the properties of
$\alpha$-pieces that can be found in \cite{hjorth}.

\begin{lem}\label{hist}(Hjort)
Let $V\in {\mathcal V}$, $U\in {\mathcal U}$, $x\in U$
and $\alpha$ be an ordinal, $\alpha>0$.
Then the following statements are true.\parskip0pt

$(1)$ $B_{\alpha}(x,U,V)$ is locally $V_U$-invariant,
$V_U x \subseteq B_{\alpha}(x,U,V)\subseteq U$\parskip0pt

\quad \ and
$B_{\alpha}(x,U,V)= B_{\alpha}(z,U,V)$ for every $z\in B_{\alpha}(x,U,V)$.
\parskip3pt

$(2)$ For any $z\in U$ the sets $B_{\alpha}(x, U,V),$ $B_{\alpha}(z,U,V)$
are either equal or disjoint.\parskip3pt

$(3)$ If $x\in U_n\subseteq U$, $V_m\subseteq V$
and $\beta\le \alpha$, then
$B_{\alpha}(x,U_n,V_m)\subseteq B_{\beta}(x,U,V)$.\parskip3pt

$(4)$ $hB_{\alpha}(x,U,V)=B_{\alpha}(hx,hU,V^h)$,
for all $h\in G_0$.\parskip3pt

\end{lem}

{\bf Remark.} While discussing $\alpha$-pieces
we may omit conditions $U_n\subseteq U$ and $V_n\subseteq V$
in the formula defining $\phi_{\alpha+1}(x,U,V)$ and let $U_n, V_n$ vary
over all elements of ${\mathcal U}$ and ${\mathcal V}$ respectively.
This is because the set $B_{\alpha+1}(x,U,V)$ coincides with the set
$$
\Big\{y:\ \{\langle\phi_{\alpha}(y',U_n,V_m), n, m\rangle:
y'\in V_Uy\cap U_n\}=\{\langle\phi_{\alpha}(x',U_n,V_m), n, m\rangle:
x'\in V_Ux\cap U_n\}\Big\}. 
$$
To see this note that the latter set is obviously included
in $B_{\alpha+1}(x,U,V)$. 

To get the converse inclusion we proceed  as follows.
Consider any $y$ which belongs to $B_{\alpha+1}(x,U,V)$ and
any triple $(x',U_n,V_m)$ with $x'\in V_Ux\cap U_n$.
Take any $U_i\subseteq U_n\cap U$ containing $x'$
and any  $V_j\subseteq V_m\cap V$.
According to the assumption on $y$
we may find $y'\in V_Uy$ such that
$\phi_{\alpha}(x',U_i,V_j)=\phi_{\alpha}(y',U_i,V_j)$, i.e.
$y'\in B_{\alpha}(x',U_i,V_j)$.
By Lemma \ref{hist}(3) the latter implies $y'\in B_{\alpha}(x',U_n,V_m)$,
i.e. $\phi_{\alpha}(x',U_n,V_m)=\phi_{\alpha}(y',U_n,V_m)$.
This proves that  the set
$\{\langle\phi_{\alpha}(x',U_n,V_m), n, m\rangle:
x'\in V_Ux\cap U_n\}$ is contained in the set 
$\{\langle\phi_{\alpha}(y',U_n,V_m), n, m\rangle:
y'\in V_Uy\cap U_n\}$. 
The symmetric argument shows that 
$$ 
\{\langle\phi_{\alpha}(y',U_n,V_m), n, m\rangle:
y'\in V_Uy\cap U_n\}\subseteq
\{\langle\phi_{\alpha}(x',U_n,V_m), n, m\rangle:
x'\in V_Ux\cap U_n\}.
$$
\bigskip

By Proposition 2.C.2 of the paper of Becker
\cite{becker} there exists a unique
partition of $X$,\ $X=\bigcup\{ Y_{t}: t\in T\}$, into invariant
$G_{\delta}$ sets $Y_{t}$ such that every $G$-orbit of $Y_{t}$ is
dense in $Y_{t}$.
To construct this partition we define for any $t\in 2^{\mathbb{N}}$ the set
$$
Y_{t}=(\bigcap\{ GA_{j}:t(j)=1\})\cap
(\bigcap\{ X\setminus GA_{j}:t(j)=0\})
$$
and take $T=\{ t\in 2^{\mathbb{N}}:Y_{t}\not=\emptyset\}$.\parskip0pt

Observe that the family $\{B_1(x,X,G):x\in X\}$
is just the canonical partition defined by Becker.
Moreover for every ordinal $0<\alpha<\omega_1$ the family
$\{B_{\alpha}(x,X,G):x\in X\}$ is a partition of $X$  approximating
the original orbit partition.
Below we will see that every such a partition also can be obtained
in a canonical way mimicking the construction of the canonical partition.

\begin{prop}\label{ph}
For every  $U\in {\mathcal U}$, $x\in U$ and $V\in {\mathcal V}$
the following equalities hold: 
$$
\begin{array}{l@{\ \ }l}
B_1(x,U,V) & =\bigcap\limits_{n}\{ V_U U_n:
V_U x\cap U_n\not=\emptyset\}\cap
\bigcap\limits_{n}\{X\setminus V_U U_n:\ V_Ux\cap U_n=\emptyset\},\\
B_{\alpha +1}(x,U,V)&=
\bigcap\limits_{n,m}\{V_UB_{\alpha}(y,U_n,V_m):
y\in U_n, V_Ux\cap B_{\alpha}(y,U_n,V_m)\not=\emptyset\}\cap\\
&\cap
\bigcap\limits_{n,m}\{X\setminus V_UB_{\alpha}(y,U_n,V_m):
y\in U_n,
V_Ux\cap B_{\alpha}(y,U_n,V_m)=\emptyset\},\\
B_{\lambda}(x,U,V) &=\bigcap
\{B_{\alpha}(x,U,V):\alpha<\lambda\}, \mbox{ for }\lambda \mbox{ limit } .
\end{array}
$$
\end{prop}
\bigskip

{\em Proof}. The first and the last equalities are obvious.
We have to prove the second one.\parskip0pt

($\subseteq$) Take any $z\in B_{\alpha +1}(x,U,V)$
and a triple $(y,U_n,V_m)$ such that  $y\in U_n$.\parskip0pt

If $V_Ux\cap B_{\alpha}(y,U_n,V_m)\not=\emptyset$,
then for some $x'\in V_Ux$ we have
$x'\in B_{\alpha}(y,U_n,V_m)$, i.e.
$\phi_{\alpha}(x',U_n,V_m)=\phi_{\alpha}(y,U_n,V_m)$.
By the assumption on $z$ there is $z'\in V_Uz$ such that
$\phi_{\alpha}(z',U_n,V_m)=\phi_{\alpha}(y,U_n,V_m)$,
i.e. $z'\in B_{\alpha}(y,U_n,V_m)$.
Hence by the properties of local saturation
$z\in V_UB_{\alpha}(y,U_n,V_m)$.\parskip0pt

On the other hand if  $V_Ux\cap B_{\alpha}(y,U_n,V_m)=\emptyset$,
then for every $x'\in V_Ux$ we have $x'\not\in B_{\alpha}(y,U_n,V_m)$,
i.e. $\phi_{\alpha}(x',U_n,V_m)\not=\phi_{\alpha}(y,U_n,V_m)$.
This implies
$\langle\phi_{\alpha}(y,U_n,V_m),n,m\rangle\not\in\phi_{\alpha+1}(x,U,V)$. 
Suppose towards contradiction that  $z\in V_UB_{\alpha}(y,U_n,V_m)$.
Then there is $z'\in V_Uz$ such that $z'\in  B_{\alpha}(y,U_n,V_m)$,
i.e. $\phi_{\alpha}(z',U_n,V_m)=\phi_{\alpha}(y,U_n,V_m)$.
Hence we get
$$
\langle\phi_{\alpha}(y,U_n,V_m),n,m\rangle\in\phi_{\alpha+1}(z,U,V). 
$$ 
Therefore we see that
$\phi_{\alpha+1}(z,U,V)\not=\phi_{\alpha+1}(x,U,V)$, thus 
$z\not\in B_{\alpha+1}(x,U,V)$. 
This contradicts our assumptions. 

($\supseteq$) Suppose that $z\not\in B_{\alpha+1}(x,U,V)$.
Then there is a pair $(n,m)$ such that
$$
\{\phi_{\alpha}(x',U_n,V_m):x'\in V_Ux\cap U_n\}\not=
\{\phi_{\alpha}(z',U_n,V_m):z'\in V_Uz\cap U_n\}. 
$$
Therefore one of the following cases holds: 

1$^o$ There is $x'\in V_Ux$ such that
$V_Uz\cap B_{\alpha}(x',U_n,V_m)=\emptyset$; 
 
2$^o$ There is $z'\in V_Uz$ such that
$V_Ux\cap B_{\alpha}(z',U_n,V_m)=\emptyset$.\\
Either case implies

$\begin{array}{l}z\not\in
\bigcap\{V_UB_{\alpha}(y,U_n,V_m):y\in U_n\subseteq U, V_m\subseteq V,
V_Ux\cap B_{\alpha}(y,U_n,V_m)\not=\emptyset\}\\
\cap\bigcap\{X\setminus V_UB_{\alpha}(y,U_n,V_m):
y\in U_n\subseteq U, V_m\subseteq  V,
V_Ux\cap B_{\alpha}(y,U_n,V_m)=\emptyset\}. \ \Box 
\end{array}$

\bigskip

From now on we shall use Proposition \ref{ph} as a definition of
an $\alpha$-piece.
\bigskip 

{\bf Remark.} Putting Proposition \ref{ph} together with Lemma \ref{hist}
we obtain the following expression
$$ 
\begin{array}{c}
B_{\alpha}(x,U,V)=\\
\bigcap\limits_{n,m\atop{\beta<\alpha}}\{V_UB_{\beta}(y,U_n,V_m):\
\ y\in U_n,
\ V_Ux\cap B_{\beta}(y,U_n,V_m)\not=\emptyset\}\cap\\
\bigcap\limits_{n,m\atop{\beta<\alpha}}\{X\setminus V_UB_{\beta}(y,U_n,V_m):\
y\in U_n,\ V_Ux\cap B_{\alpha}(y,U_n,V_m)=\emptyset\}.
\end{array}
$$
Hence we may formulate the following assertion.

\begin{cor}\label{phar} 
Let $x,y\in U$.
Then $B_{\alpha}(x,U,V)=B_{\alpha}(y,U,V)$ if and only if
the local orbits $V_Ux$ and $V_Uy$
intersect exactly the same $\beta$-pieces,
for every $0\le \beta<\alpha$.\parskip3pt
\end{cor} 

In the next lemma we formulate another important property of
$\alpha$-pieces: every element of $V_UB_{\alpha}(x, U_n, V_m)$
can be surrounded by some $\alpha$-piece entirely contained
in $V_UB_{\alpha}(x, U_n, V_m)$.
This property will be frequently  applied in the rest of the paper.

\begin{lem}\label{list}
Let $V\in {\mathcal V}$, $U\in {\mathcal U}$ and $x\in U_n\cap U$.
Then for every $m\in\omega$, ordinal $\alpha > 0$
and $y\in V_UB_{\alpha}(x,U_n,V_m)$
there are $U_i\subseteq U_n$ and $h\in G_0$ such that
$h\in \langle V\rangle ^{U_i}_U$ (i.e.
$h\in {\langle V\rangle}^t_U$ for every $t\in U_i$),
$y\in hU_i$ and\quad
$B_{\alpha}(y,hU_i, (V_m)^h)\subseteq V_UB_{\alpha}(x,U_n,V_m)$.
\end{lem}

{\em Proof.} If $y\in V_UB_{\alpha}(x,U_n,V_m)$,
then the local orbit $V_Uy$ intersects $B_{\alpha}(x,U_n,V_m)$.
Thus the intersection
$\langle V\rangle_U^y\cap \{g\in G: gy\in B_{\alpha}(x,U_n,V_m)\}$
is nonempty.
By Lemma \ref{bH}(3) we see that  $\langle V\rangle_U^y$ is open.
The set $\{g\in G: gy\in B_{\alpha}(x,U_n,V_m)\}$
is open either.
Indeed if  $gy\in B_{\alpha}(x,U_n,V_m)$,
then by Lemma \ref{hist} the local orbit
$(V_m)_{U_n}gy$  is entirely contained in $B_{\alpha}(x,U_n,V_m)$
and so
$\langle V_m\rangle_{U_n}^{gy}g\subseteq
\{g\in G: gy\in B_{\alpha}(x,U_n,V_m)\}$.
Thus we are done since the set $\langle V_m\rangle_{U_n}^{gy}g$
is open by Lemma \ref{bH}(3). 

To sum up
$\langle V\rangle_U^y\cap \{g\in G: gy\in B_{\alpha}(x,U_n,V_m)\}$
is a nonempty open set.
Therefore it contains an element of $G_0$,
i.e. we can find $\ g\in G_0\cap \langle V\rangle_U^y\ $ such that
$gy\in B_{\alpha}(x,U_n,V_m)$.
Put $y'=gy$ and $h=g^{-1}$.
We have $h\in G_0$ and  $hy'=y$.
By Lemma \ref{bH}(2) we see that $h\in \langle V\rangle_U^{y'}$.
According to Lemma \ref{bH}(3) there is a basic open $U_i\subseteq U_n$
containing $y'$ such that  $h\in {\langle V\rangle}^{U_i}_U$.
Since $B_{\alpha}(y',U_i,V_m)\subseteq U_i$ and
$B_{\alpha}(y',U_i,V_m)\subseteq B_{\alpha}(y',U_n,V_m)=
B_{\alpha}(x,U_n,V_m)$,
then 
$$
B_{\alpha}(hy',hU_i,(V_m)^h)=hB_{\alpha}(y',U_i,V_m)
\subseteq V_UB_{\alpha}(x,U_n,V_m). 
$$
$\Box$
\bigskip

Now we are ready to formulate the main result of this part.
Despite its technical character
this theorem shed a new light on the nature of
$\alpha$-pieces.
In particular it enables us to prove that $\alpha$-pieces
are Borel sets.

\begin{thm}\label{VB}  Let $U\in {\mathcal U}$, $V\in {\mathcal V}$
and $x\in U$.
Then for every ordinal $\alpha>0$ the following equality is true
$$ 
\begin{array}{c}
B_{\alpha+1}(x,U,V)=\\
\bigcap\limits_{n,m}
(\bigcup\{B_{\alpha}(gx,hU_i,(V_m)^h): U_i\subseteq U_n,
h\in  G_0\cap {\langle V\rangle}^{U_i}_U,
g\in G_0\cap {\langle V\rangle}^x_U, gx\in hU_i\})\\
\cap\bigcap\limits_{n,m}\big(
(X\setminus U_n)\cup
\bigcup\{B_{\alpha}(gx,U_n,V_m):
g\in G_0\cap {\langle V\rangle}^x_U, gx\in U_n\}\big).
\end{array}
$$ 
\end{thm}

{\em Proof}. We will apply Proposition \ref{ph}.
Consider  arbitrary $U_n$, $V_m$ and $y\in U_n$
such that $B_{\alpha}(y,U_n,V_m)\cap V_Ux\not=\emptyset$.
The set $V_UB_{\alpha}(y,U_n,V_m)$ can be presented as a union
of $\alpha$-pieces in the way described in Lemma \ref{list}.
If we throw aside the elements that can be surrounded by some $\alpha$-piece
disjoint from $V_Ux$, then we may limit ourselves to $\alpha$-pieces
containing elements of the form $gx$,
where $g\in G_0\cap {\langle V\rangle}^x_U$.
As a result we obtain the following formula.
$$\begin{array}{c}
V_UB_{\alpha}(y,U_n,V_m)\cap
\bigcap\limits_{i,j} \{(X\setminus V_UB_{\alpha}(z,U_i,V_j)):
V_Ux\cap B_{\alpha}(z,U_i,V_j)=\emptyset\} \\
=\\
\bigcup\{B_{\alpha}(gx,hU_i,(V_m)^h): U_i\subseteq U_n,
h\in  G_0\cap {\langle V\rangle}^{U_i}_U,
g\in G_0\cap {\langle V\rangle}^x_U, gx\in hU_i\}\cap\\
\bigcap\limits_{i,j} \{(X\setminus V_UB_{\alpha}(z,U_i,V_j)):
V_Ux\cap B_{\alpha}(z,U_i,V_j)=\emptyset\}.
\end{array}
$$
Next consider the intersection
$$
\bigcap\limits_{n,m}\{X\setminus V_UB_{\alpha}(y,U_n,V_m):
\ y\in U_n, V_Ux\cap B_{\alpha}(y,U_n,V_m)=\emptyset\}. 
$$
By Lemma \ref{list} it is a complement of a union
$$
\bigcup\{B_{\alpha}(y,U_n,V_m): y\in U_n, V_Ux\cap B_{\alpha}(y,U_n,V_m)
=\emptyset\}
$$
To complete the proof observe that by Lemma \ref{hist}(2)
for any  $n,m$ and $y\in U_n$ the following equivalence is true:
$$
B_{\alpha}(y,U_n,V_m)\cap V_Ux=\emptyset \mbox{ \  if and only if}
$$ 
$$
B_{\alpha}(y,U_n,V_m)\subseteq U_n\setminus
\bigcup\{B_{\alpha}(gx,U_n,V_m):
g\in G_0\cap {\langle V\rangle}^x_U, gx\in U_n\}. 
$$
$\Box$
\bigskip

Involving Theorem \ref{VB} in an inductive argument we can prove
the following statement.

\begin{cor}\label{borel}
Let $x,y\in X$, $V\in {\mathcal V}$, $U\in {\mathcal U}$,
\  $\ \alpha\ge 1\ $ be an ordinal and
$$
\rho(\alpha)=\left\{\begin{array}{l@{\quad\mbox{if}\quad}l}
2n  & \alpha=n, \mbox{ where }n<\omega\\
\beta+2n & \alpha=\beta+n, \mbox{ where }n<\omega
\mbox{ and }\beta\mbox{ is limit}.
\end{array}\right.
$$
Then we have
$B_{\alpha}(x,U,V)\in \Pi^0_{\rho(\alpha)}(X)$.
\end{cor}
\bigskip

Finally it follows from Lemma \ref{hist}(2) and  Theorem \ref{VB}
that for every ordinal $\alpha\ge 1$
the family of $\{B_{\alpha}(x,X,G):x\in X\}$
forms a partition of the space $X$ into invariant Borel sets
which Borel rank is bounded by a countable ordinal.

\begin{cor} \label{complex}
For every $V\in {\mathcal V}$ and $U\in {\mathcal U}$
the family
$\{B_{\alpha}(x,U,V):x\in U\}$ is a partition of $U$ into
locally $V_U$-invariant $\Pi^0_{\rho(\alpha)}$-sets.
In particular the family
$\{B_{\alpha}(x,X,G):x\in X\}$ is a partition of the whole space $X$ into
$G$-invariant $\Pi^0_{\rho(\alpha)}$-sets.
\end{cor}

In the end of this section we shall prove  another property
of $\alpha$-pieces. 
Lemma \ref{hist}(4) states that
for any  $\alpha$-piece $B_{\alpha}(x,U,V)$ and any $h\in G_0$,
the set $hB_{\alpha}(x,U,V)$ is an $\alpha$-piece defined with respect
to basic open $hU,V^h$.
Since $\alpha$-pieces are defined only with respect
to basic open $U,V$,  the above is not true
in the general case of any $h\in G$.
Instead we can prove a related property which can be viewed as a
generalization of Lemma \ref{hist}(3).

\begin{lem} Let $\alpha>0$ be an ordinal, $U,U'\in {\mathcal U}$,
$V,V'\in {\mathcal V}$, $x\in U$ and $h\in G$.
If $hx\in U'\subseteq hU$ and $V'\subseteq V^h$,
then  for every ordinal $\beta\le\alpha$ we have
$$B_{\alpha}(hx,U',V')\subseteq hB_{\beta}(x,U,V).$$
\end{lem}

{\em Proof}. First observe that according to Lemma \ref{hist}(3)
we have to consider only the case $\alpha=\beta$.
We proceed by induction on $\alpha$ applying Proposition \ref{ph}.
\parskip3pt

Assume $\alpha=1$.  Suppose towards contradiction that
$B_1(hx,U',V')\not\subseteq hB_1(x,U,V)$.
Then there is $y\in B_1(hx,U',V')$ which does not belong to
$hB_1(x,U,V)$.
Hence by the assumption that $U'\subseteq hU$ we have
$h^{-1}y\in U$ and $h^{-1}y\not\in B_1(x,U,V)$.
So according to Lemma \ref{hist}(2) we see that
$B_1(x,U,V)\cap B_1(h^{-1}y,U,V)=\emptyset$.
Then there is a basic open set $U_n$ meeting one of the local orbits
$V_Ux$, $V_U(h^{-1}y)$ and disjoint from the other.
W.l.o.g. we may assume that
$V_Ux\cap U_n\not=\emptyset$ while $V_U(h^{-1}y)\cap U_n=\emptyset$.
This implies $V_Ux\subseteq V_UU_n$
and  $V_U(h^{-1}y)\cap V_UU_n=\emptyset$.
Using Lemma \ref{locsat}(6) we get
$(V^h)_{hU}hx\subseteq (V^h)_{hU}hU_n$ while
$(V^h)_{hU}y\cap (V^h)_{hU}hU_n=\emptyset$.
Then by Lemma \ref{locsat}(1) we have
$(V')_{U'}hx\subseteq (V^h)_{hU}hU_n$ and
$(V')_{U'}y\cap (V^h)_{hU}hU_n=\emptyset$.
Since according to Lemma \ref{bH}(3) the set $(V^h)_{hU}hU_n$ is open,
the latter conditions imply
that it contains a basic open set meeting $(V')_{U'}hx$
and disjoint from $(V')_{U'}y$.
This contradicts the assumption that $y\in B_1(x,U',V')$.
\parskip3pt

For the successor step assume that  for every $n,m,k,l\in\omega$
and $h\in G$  satisfying
the conditions $hx\in U_n\subseteq hU_k$ and $V_m\subseteq V_l^h$
we have
$B_{\alpha}(hx,U_n,V_m)\subseteq hB_{\alpha}(x,U_k,V_l)$.
Suppose towards contradiction that
$B_{\alpha+1}(hx,U',V')\not\subseteq hB_{\alpha+1}(x,U,V)$.
There is $y\in B_{\alpha+1}(hx,U',V')$ such that
$y\not\in hB_{\alpha+1}(x,U,V)$.
Then $h^{-1}y\not\in B_{\alpha+1}(x,U,V)$ and so
$B_{\alpha+1}(x,U,V)$ is disjoint from $B_{\alpha+1}(h^{-1}y,U,V)$.
By Proposition \ref{ph} there is an $\alpha$-piece
$B_{\alpha}(z,U_k,V_l)$ meeting exactly one of the local orbits
$V_Ux$ and $V_U(h^{-1}y)$.
W.l.o.g. assume that $B_{\alpha}(z,U_k,V_l)$ meets $V_Ux$.
We have $V_Ux\subseteq V_UB_{\alpha}(z,U_k,V_l)$
and  $V_U(h^{-1}y)\cap V_UB_{\alpha}(z,U_k,V_l)=\emptyset$.
Since $x\in V_UB_{\alpha}(z,U_k,V_l)$, then according to Lemma \ref{list}
there are $i,j\in \omega$ such that
$B_{\alpha}(x,U_{i},V_{j})\subseteq V_UB_{\alpha}(z,U_k,V_l)$.
Hence $B_{\alpha}(x,U_{i},V_{j})$ is disjoint from $V_U(h^{-1}y)$ either,
and so we may assume that $x=z$, $i=k$ and $j=l$.\parskip0pt

Now take some $U_n, V_m$ so that $hx\in U_n\subseteq hU_k$
and $V_m\subseteq V_l^h$.
We claim that the local orbit $(V')_{U'}y$
is disjoint from $B_{\alpha}(hx,U_n,V_m)$.
Otherwise the local orbit $(V^h)_{hU}y$ is not disjoint from
$B_{\alpha}(hx,U_n,V_m)$ either.
By the inductive assumption we have
$B_{\alpha}(hx,U_n,V_m)\subseteq hB_{\alpha}(x,U_k,V_l)$.
Hence $(V^h)_{hU}y$ meets $hB_{\alpha}(x,U_k,V_l)$
which by Lemma \ref{locsat}(6) implies that
$V_U(h^{-1}y)$ meets  $B_{\alpha}(x,U_k,V_l)$.
This contradiction shows that the claim that
$(V')_{U'}y\cap B_{\alpha}(hx,U_n,V_m)=\emptyset$ is true.
\parskip0pt

Hence we can apply Proposition \ref{ph} to conclude that
$y\not\in  B_{\alpha+1}(hx,U',V')$.
This contradicts our assumption and completes the successor step. 
The limit step is immediate. $\ \Box$ 
\bigskip

Using the lemma above we can prove the following assertion.
It will be applied in the proof of Theorem \ref{pol}.

\begin{cor}\label{basis}
For every $U\in {\mathcal U}$, $V\in {\mathcal V}$, $x\in U$,
$h\in G$ and ordinal $\alpha>0$ the set
$hB_{\alpha}(x,U,V)$ is a union of appropriate $\alpha$-pieces.
\end{cor}

{\em Proof}. Take any $y\in hB_{\alpha}(x,U,V)$.
According to Lemma \ref{hist}(2) we may assume that $y=hx$.
For some $n,m\in\omega$  we have $hx\in U_n\subseteq hU$ and
$(V_m)\subseteq V^h$.
Then by the lemma above we get
$B_{\alpha}(hx,U_n,V_m)\subseteq hB_{\alpha}(x,U,V)$.
$\Box$
\bigskip 

The results above starting from Proposition \ref{ph} till
Theorem \ref{VB} and its consequences concerning borelness of
the sets $B_{\alpha}(x,U,V)$ is our contribution
in study of $\alpha$-pieces.

\subsection{Finer topologies}

Since every piece of the canonical partition
is a $G_{\delta}$-subset of $X$,
it is a Polish space with the topology inhertited from the original
Polish topology on $X$.
This fact is generalized by Hjorth (see \cite{hjorth}) who proved
that for every $x\in X$ the set $B_{\gamma^{\star}(x)+2)}(x,X,G)$
is a Polish $G$-space with respect to the topology generated
by the family
$\{B_{\gamma^{\star}(x)}(x',U_n,V_m): x'\in Gx\cap U_n,\ n,m\in\omega\}$.
We improve this result and show that for every ordinal $\alpha>0$
every $\alpha$-piece of the form $B_{\alpha}(x,X,G)$
is a Polish $G$-space with respect to the topology generated by 'ealier'
$\beta$-pieces  and this topology generates
the same Borel structure as the original topology.
Our proof is based on Theorem \ref{VB} and the theorem on Borel families
by Sami. 

From now on we shall use the following notation for every ordinal $\beta$:
$$
\begin{array}{l@{\ =\ }l}
{\mathcal B}^x_0&{\mathcal U}, \\ 
{\mathcal B}^x_{\beta}& \{B_{\beta}(x',U,V): U\in {\mathcal U},
x'\in U\cap G_0x, V\in {\mathcal V}\}
\mbox{\quad for\quad }\beta>0, \\
{\mathcal B}^x_{<\beta}&\bigcup\{{\mathcal B}^x_{\gamma}:\ \gamma<\beta \} .  
\end{array}
$$

\begin{thm}\label{pol}
Let $x\in X$ and $0<\alpha<\omega_1$ be an ordinal.
The set $B_{\alpha}(x,X,G)$
with the (relative) topology generated by the family
${\mathcal B}^x_{<\alpha}$ as basic open sets is a Polish $G$-space
with the same Borel structure as the original topology.\parskip0pt

Additionally for every $U\in {\mathcal U}$,  $V\in {\mathcal V}$
with $x\in U$ the set $B_{\alpha}(x,U,V)$ is a Polish space with respect to
this topology.
\end{thm}

{\em Proof.}
As we have already mentioned $B_1(x,X,G)$ is a $G_{\delta}$ subset
of $\langle X, \tau\rangle $,
thus it is a Polish space with respect to the (relative) topology
generated by ${\mathcal U}$.
Therefore below we will deal only with $\alpha>1$.
We shall use the following result by Sami (see \cite{sami}, Lemma 4.2).
\begin{quote}
{\em Let $\langle X, t\rangle$ be a topological space and
$0\le \zeta<\omega_1$.
Let ${\mathcal F}$ be a {\em Borel family of rank} $\zeta$, i.e.
a family of subsets of $X$ which can be decomposed into
subfamilies of two types
$\ {\mathcal F}=\bigcup\{P_{\xi}:0\le\xi<\zeta\}\cup
\bigcup\{S_{\xi}:0\le\xi<\zeta\}
\ $
satisfying the following conditions: \\ 
1. $S_0$ consists of open sets,\\ 
2. $P_{\xi}=\{X\setminus A:\ A\in S_{\xi}\}$, for $0\le \xi<\zeta$,\\ 
3. every element of $S_{\xi}$ is a union of a countable subfamily
of $\bigcup \{P_{\eta}:0\le\eta<\xi\}$, for $0\le \xi<\zeta$.\parskip0pt

If $\langle X, t\rangle$ is a Polish space
then the topology generated by a family
of intersections of finite subsets of the union $t\cup {\mathcal F}$ is also
Polish.}
\end{quote}

We start with some preliminary work.
For every $0\le \xi<\omega_1$ we define the sets $S_{\xi}$ and $P_{\xi}$.
First we put:
$$
\begin{array}{r@{\ }l}
S_0=&\{V_U A: A, U\in {\mathcal U}, V\in {\mathcal V}\}, \\
P_0=&\{X\setminus D:\  D\in S_0\}, \\
& \\
S_1=&P_0 , \\
P_1=&\{X\setminus D: \ D\in S_1\}=S_0 , \\
\\
S_2=&\{(\bigcup\{X\setminus V_{U} A:
\ A\in {\mathcal U},\ x'\in V_{U}A\}):
\ U\in {\mathcal U},\ V\in {\mathcal V}, x'\in G_0x\}\\
&\quad \cup \, \{(\bigcup\{V_U A:
\ A\in {\mathcal U},\  x'\not\in {V}_{U}A\}):
 \ U\in {\mathcal U},\ V\in {\mathcal V}, x'\in G_0x\} , \\
P_2=&\{X\setminus D:\ D\in S_2\}. 
\end{array}
$$  
Observe that
$\bigcup\{(S_i\cup P_i): 0\le i < 2\}$
is a Borel family of rank 3, where 
$$
\begin{array}{l}
S_2=\Big\{(X\setminus B_1(x',U,V)):
\ U\in {\mathcal U}, x'\in U, \ V\in {\mathcal V}, x'\in G_0x\Big\} , \\
P_2={\mathcal B}_1^x.
\end{array}
$$ 
We proceed similarly at each successor stage.
Every  successor ordinal has one of
the forms:   $\xi + 2n+1$ or $\xi +2n+2$,
where $n$ is a natural number and
$\xi=0$ or $\xi$ is a limit ordinal.
We define: 
$$
\begin{array}{r} S_{\xi +2n+1}=
\Big\{\Big(\bigcup\{B_{\xi+n}(gx',hU_i,(V_m)^h):\ U_i\subseteq U_n,
\ h\in G_0\cap {\langle V\rangle}_U^{U_i},
g\in G_0\cap {\langle V\rangle}_U^{x'},\\
gx'\in hU_i \}\Big):
\ n,m\in\omega, U\in {\mathcal U}, V\in {\mathcal V}, x'\in G_0x\Big\}
\\
\cup\, \Big\{\Big((X\setminus U_n)\cup
\bigcup\{  B_{\xi+n}(gx',U_n,V_m): g\in G_0\cap {\langle V\rangle}_U^{x'},
gx'\in U_n\}\Big):\\
n,m\in\omega, U\in {\mathcal U}, V\in {\mathcal V}, x'\in G_0x\Big\} , \\
\end{array}
$$
$$
\begin{array}{l}
P_{\xi +2n+1}= \{X\setminus D:D\in S_{\xi+2n+1}\} , \\
\\
S_{\xi +2n+2}=\{(X\setminus B): B\in {\mathcal B}^x_{\xi +n+1}\} , \\
P_{\xi +2n+2}={\mathcal B}^x_{\xi +n+1} .
\end{array}
$$ 
Finally, for every limit $\xi<\omega_1$ we put 
$$
\begin{array}{l}
S_{\xi}=\{X\setminus B:\ B\in {\mathcal B}^x_{\xi}\} ,\\
P_{\xi}={\mathcal B}^x_{\xi}.
\end{array}
$$ 
We claim that for every $1\le \zeta<\omega_1$ the family
$\bigcup\{(S_{\xi}\cup P_{\xi})\ :\xi<\zeta\}$ is a Borel family
of rank $\zeta$.
It is clear that such a family satisfies conditions 1-2 of the Sami's theorem.
We have to check that  it also satisfies condition 3.
We apply an inductive argument.
It is obvious for $\zeta=1, 2$.
The case of a limit $\zeta$ immediately follows from the inductive assumption.

For the successor step take an arbitrary $\zeta\le 1$ and
suppose that the family
$\bigcup\{(S_{\xi}\cup P_{\xi}):\ 0\le\xi<\zeta\}$
satisfies condition 3.
We have to check that the family
$\bigcup\{(S_{\xi}\cup P_{\xi}): 0\le\xi<\zeta+1\}$
also satisfies condition 3.
Since 
$$
\bigcup\{(S_{\xi}\cup P_{\xi}): 0\le\xi<\zeta+1\}=
\bigcup\{(S_{\xi}\cup P_{\xi}): 0\le \xi<\zeta\}
\ \cup \ S_{\zeta}\cup P_{\zeta}, 
$$
it suffices to prove that every element from
$S_{\zeta}$ is a union of elements of the set
$\bigcup \{P_{\xi}: 0\le\xi<\zeta\}$.
Consider two cases.\parskip3pt

$1^o$ $\zeta$ is a successor ordinal.
Then there are unique ordinals $\gamma$ and $n$ such that
$n$ is a  natural number,
$\gamma$ equals $0$ or is a limit ordinal and
$\zeta$  has one of the following form:
$\gamma+2n+1$ or $\gamma+2n+2$.
In the first case the desired property follows by the definition.
For the second case consider any element
$D=X\setminus B_{\gamma+n+1}(x',U,V)$ from $S_{\gamma+2n+2}$.
Applying Theorem \ref{VB} we have
$$
\begin{array}{c}D=\bigcup\limits_{n,m}
(X\setminus\bigcup\{B_{\alpha}(gx',hU_i,(V_m)^h): U_i\subseteq U_n,
h\in  G_0\cap {\langle V\rangle}^{U_i}_U,
g\in G_0\cap {\langle V\rangle}^{x'}_U, gx'\in hU_i\})\\
\cup\bigcup\limits_{n,m}\Big(X\setminus
((X\setminus U_n)\cup
\bigcup\{B_{\alpha}(gx',U_n,V_m):
g\in G_0\cap {\langle V\rangle}^{x'}_U, gx'\in U_n\})\Big).
\end{array}
$$
Hence we see that $D$ is a union of elements from $P_{\gamma+2n+1}$. 

$2^o$ $\zeta$ is limit.
Then by  Proposition \ref{ph} we have 
$$
\begin{array}{l} 
P_{\zeta}=
\{B:\ B\in{\mathcal B}^x_{\zeta}\}=
\{ \bigcap\limits_{\xi<\zeta} B_{\xi}(gx,U,V):
g\in G_0, U\in {\mathcal U}, V\in {\mathcal V}\} , \\
S_{\zeta}=\{X\setminus B:\ B\in {\mathcal B}^x_{\zeta}\}=
\{ \bigcup\limits_{\xi<\zeta}
(X\setminus B_{\xi}(gx,U,V)):
g\in G_0, U\in {\mathcal U}, V\in {\mathcal V}\}.
\end{array}
$$

To show that every element of $S_{\zeta}$ is a union of elements
from $\bigcup\{P_{\xi}:\xi<\zeta\}$ we shall use the following property.
\parskip3pt

{\em Claim 1}. If $\xi$ equals $0$ or is a limit ordinal
and $n$ is a natural number then
${\mathcal B}^x_{\xi+n}\subseteq S_{\xi+2n+1}$.
\parskip2pt

{\em Proof of Claim 1}. Consider an arbitrary element $B_{\xi+n}(x',U,V)$
from ${\mathcal B}^x_{\xi+n}$.
If we subsitute $U_n=U$ and $V_m=V$ in the first formula defining
elements of $S_{\xi+2n+1}$, then we get the set
$$
D=\bigcup\{B_{\xi+n}(gx',hU_i,(V)^h):U_i\subseteq U
\ h\in G_0\cap {\langle V\rangle}_U^{U_i},
g\in G_0\cap {\langle V\rangle}_U^{x'}, gx'\in hU_i \}). 
$$
The condition  
$(\star)\ U_i\subseteq U\ \wedge
\ h\in G_0\cap {\langle V\rangle}_U^{U_i}\ \wedge
\ g\in G_0\cap {\langle V\rangle}_U^{x'}\ \wedge \  gx'\in hU_i $
and Lemma \ref{hist} (1)-(3) imply
that  $D\subseteq B_{\xi+n}(x',U,V)$.
On the other hand since $h=g=1_G$ and $U_i=U$ also satisfy
condition $(\star)$, thus we get  $D\supseteq B_{\xi+n}(x',U,V)$. 

Applying this property  and the assumption that
$\zeta$ is limit we see that
${\mathcal B}^x_{<\zeta}\subseteq \bigcup\limits_{\xi<\zeta}S_{\xi}$.
Hence $\{X\setminus B:\  B\in {\mathcal B}^x_{<\zeta}\}
\subseteq\bigcup\limits_{\xi<\zeta}P_{\xi}$.
Therefore every element of $S_{\zeta}$ is a countable union of elements
of the set $\bigcup\limits_{\xi<\zeta}P_{\xi}$. 

This completes Case $2^o$. 

Now let $\gamma$ and $k$ be the unique ordinals such that
$\alpha=\gamma+k$, $k$ is a natural number and
$\gamma$ equals $0$ or is a limit ordinal.
Define
$$
\hat{\alpha}=
\left\{\begin{array}{l@{\ \mbox{if}\ }l}
\gamma+2(k-1)+1 & k>0\\
\alpha & k=0.\end{array}\right.
$$
Put 
${\mathcal F}_{\alpha}=\bigcup\{P_{\xi}\cup
S_{\xi}:\ 0\le\xi<\hat{\alpha}\}$.
We have proved that ${\mathcal F}_{\alpha}$
is a Borel family of  rank $\hat{\alpha}$ and
${\mathcal U}\subseteq {\mathcal F}_{\alpha}$.
Hence it is a subbase of the Polish topology finer then
the initial topology generated by ${\mathcal U}$.
Since ${\mathcal B}^x_{<\alpha}\subseteq {\mathcal F}_{\alpha}$,
then $B_{\alpha}(x,X,G)$
is a $G_{\delta}$ set with respect to this
topology.
Thus every $B_{\alpha}(x,X,G)$ is also a Polish space with
the inherited topology.
We now show that the family
$\{B_{\alpha}(x,X,G)\cap B:\ B\in {\mathcal B}^x_{<\alpha}\}$
is a basis of the topology.
So we have to prove that every set of the form
$B_{\alpha}(x,X,G)\cap D$, where $D\in {\mathcal F}_{\alpha}$
is a union of elements from
$\{B_{\alpha}(x,X,G)\cap B:\ B\in {\mathcal B}^x_{<\alpha}\}$.
This is an immediate consequence of the following claim.\parskip3pt

{\em Claim 2.} Let $\zeta<\beta<\alpha$.
Then for every $x'\in G_0x$, $U\in {\mathcal U}$, $V\in {\mathcal V}$
and $n,m\in\omega$ the sets below are unions of elements from the family
$\{B_{\alpha}(x,X,G)\cap B:\ B\in  {\mathcal B}^x_{\beta}\}$:  
$$
\begin{array}{l}
B_{\alpha}(x,X,G)\cap B_{\zeta}(x',U_n,V_m),\\
\\
B_{\alpha}(x,X,G)\setminus
(\bigcup\{B_{\xi+n}(gx',hU_i,(V_m)^h):\\
U_i\subseteq U_n,
\ h\in G_0\cap {\langle V\rangle}_U^{U_i},
g\in G_0\cap {\langle V\rangle}_U^{x'}, gx'\in hU_i \}),\\
\\
B_{\alpha}(x,X,G)\cap (U_n\setminus
\bigcup\{  B_{\zeta}(hx',U_n,V_m): h\in G_0\cap {\langle V\rangle}_U^{x'}\}). 
\end{array}
$$

{\em Proof of Claim 2.}
Take any $y\in B_{\alpha}(x,X,G)\cap B_{\zeta}(x',U_n,V_m)$.
By Lemma \ref{hist}(1), (3) we see that
$B_{\alpha}(x,X,G)\cap B_{\beta}(y,U_n,V_m)\subseteq
B_{\alpha}(x,X,G)\cap B_{\zeta}(x',U_n,V_m)$.
Since $y\in B_{\alpha}(x,X,G)$, then we may apply Corollary \ref{phar}
to find some $x''\in G_0x$ such that
$B_{\beta}(x'',U_n,V_m)=B_{\beta}(y,U_n,V_m)$.
\parskip0pt

To settle the second part of this claim consider any
$y\in B_{\alpha}(x,X,G)$ which does not belong to the union
$$
\begin{array}{c}\bigcup\{B_{\xi+n}(gx',hU_i,(V_m)^h):U_i\subseteq U_n,
\ h\in G_0\cap {\langle V\rangle}_U^{U_i},
g\in G_0\cap {\langle V\rangle}_U^{x'}, gx'\in hU_i \}.
\end{array}
$$
Applying the argument from the proof of Theorem \ref{VB} we see that
$$
\begin{array}{c}
y\in B_{\alpha}(x,X,G)\setminus V_UB_{\zeta}(x',U_n,V_m)\end{array}.
$$
By Lemma \ref{hist}(3) and  Proposition \ref{ph} we have
$B_{\beta}(y,U,V)\subseteq X\setminus V_UB_{\zeta}(x',U_n,V_m)$.
We now  find  $x''\in G_0x$ such that
$B_{\beta}(x'',U_n,V_m)=B_{\beta}(y,U_n,V_m)$ and finish the proof.
\parskip0pt

Similarly we prove the last part of the claim.

We now see that
$\{B_{\alpha}(x,X,G)\cap B:\ B\in  {\mathcal B}^x_{<\alpha}\}$
generates on $B_{\alpha}(x,X,G)$ the (relative) topology
defined by ${\mathcal F}_{\alpha}$.
Since every $\alpha$-piece $B_{\alpha}(x,U,V)$ is a $G_{\delta}$
subset of $B_{\alpha}(x,X,G)$ we see that the additional  statement of
this theorem is true either. 

Now it suffices to show that the action
$a: G \times B_{\alpha}(x,X,G)\to B_{\alpha}(x,X,G)$
is continuous with respect to each coordinate.
Take an arbitrary basic open set
$B=B_{\beta}(x',U,V)\cap B_{\alpha}(x,X,G)$,
where $\beta<\alpha$, $U\in {\mathcal U}$, $V\in {\mathcal V}$
and $x'\in G_0x\cap U$.
To prove continuity with respect to the first coordinate
fix some $y\in B_{\alpha}(x,X,G)$ and consider the set
$\{h\in G: hy\in B\}$.
If $h$ is an element of this set,
i.e. $hy\in  B$,
then according to Lemma \ref{hist}(1) $V_Uhy\subseteq B$.
Hence ${\langle V\rangle}_U^{hy}h$
is an open neighbourhood of $h$ contained in  $\{h\in G: hy\in B\}$.
\parskip0pt

To prove continuity with respect to the second coordinate
fix some $h\in G$ and consider the set
$$
\{y\in B_{\alpha}(x,X,G): hy\in B\}=
h^{-1}B_{\beta}(x',U,V)\cap B_{\alpha}(x,X,G).
$$
By Corollary \ref{basis} the set  $h^{-1}B_{\beta}(x',U,V)$
is a union of $\beta$-pieces.
Hence $h^{-1}B_{\beta}(x',U,V)\cap B_{\alpha}(x,X,G)$ is a union
of $\beta$-pieces meeting $V_Ux$, thus it is open with respect
to  $t^x_{\alpha}$.
This proves continuity with respect to the second coordinate.
$\Box$
\bigskip

From now on let $t_{\alpha}^x$ denote the Polish topology on
$B_{\alpha}(x,X,G)$ described above.
Observe that in the case when $\alpha$ is a successor ordinal
and $\alpha=\beta+1$, the topology $t_{\alpha}^x$ is also (relatively)
generated by a smaller basis, namely ${\mathcal B}^x_{\beta}$.
It follows directly from Claim 2.\parskip3pt

Since $t_{\alpha}^x$  is finer then the original (relative)
topology on $B_{\alpha}(x,X,G)$ all operations and sets introduced
so far can be considered with respect to $t^x_{\alpha}$.
We shall use the superscript $^{t^x_{\alpha}}$ to stress that
a given object is constructed in the $G$-space
$B_{\alpha}(x,X,G)$ with respect to the topology $t^x_{\alpha}$.
Let us illustrate this idea.\parskip3pt

{\bf Example.} Take arbitrary $x\in X$ and an ordinal $\alpha>0$.
Consider $B_{\alpha}(x,X,G)$ as a $G$-space with respect to the topology
$t^x_{\alpha}$.
Fix some enumeration $\{D_n:n\in\omega\}$ of its basis
${\mathcal U}^x_{\alpha}=
\{B_{\alpha}(x,X,G)\cap B:\ B\in  {\mathcal B}^x_{<\alpha}\}$.
Then for every $\beta>0$,
$D\in {\mathcal U}^x_{\alpha}$ and $V\in{\mathcal V}$,
we define the $\beta$-piece $B_{\beta}^{t^x_{\alpha}}(y,D,V)$ with respect
to $t^x_{\alpha}$ using the scheme from Proposition \ref{ph}
in the following way:
$$
\begin{array}{l@{\ \ }l}
B^{t^x_{\alpha}}_1(z,D,V) & =\bigcap\limits_{n}\{ V_D D_n:
V_D z\cap D_n\not=\emptyset\}\cap
\bigcap\limits_{n}\{B_{\alpha}(x,X,G)\setminus V_D D_n:\ V_Dz\cap D_n=\emptyset\},\\
B^{t^x_{\alpha}}_{\beta +1}(z,D,V)&=
\bigcap\limits_{n,m}\{V_DB^{t^x_{\alpha}}_{\beta}(y,D_n,V_m):
y\in D_n, V_Dz\cap B^{t^x_{\alpha}}_{\beta}(y,D_n,V_m)\not=\emptyset\}\cap\\
&\cap
\bigcap\limits_{n,m}\{B_{\alpha}(x,X,G)\setminus V_DB^{t^x_{\alpha}}_{\beta}(y,D_n,V_m):
y\in D_n, V_Dz\cap B^{t^x_{\alpha}}_{\beta}(y,D_n,V_m)=\emptyset\},\\
B^{t^x_{\alpha}}_{\lambda}(z,D,V) &=\bigcap
\{B^{t^x_{\alpha}}_{\beta}(z,D,V):\beta<\lambda\}, \mbox{ for }
\lambda \mbox{ limit } .
\end{array}
$$

There is a natural relationship between $\alpha$-pieces constructed
with respect to the subsequent topologies.

\begin{prop}\label{subs} Let $V\in {\mathcal V}$, $U\in {\mathcal U}$,
$x\in X$ and $x'\in Gx\cap U$.
Let $\gamma,\alpha$ be ordinals such that
$1\le \gamma<\alpha<\omega_1$.
Then for every
$y\in B_{\alpha}(x,X,G)\cap B_{\gamma}(x',U,V)$ and $\beta<\omega_1$
the following equality holds.
$$
B_{\alpha+\beta}(y,U,V)=
\left\{\begin{array}{l@{\quad\mbox{if}\quad}l}
B_{\beta+1}^{t^x_{\alpha}}(y,B_{\alpha}(x,X,G)\cap B_{\gamma}(x',U,V),V)
&\beta\mbox{ is finite }\\
B_{\beta}^{t^x_{\alpha}}(y,B_{\alpha}(x,X,G)\cap B_{\gamma}(x',U,V),V)
&\beta\mbox{ is infinite}.
\end{array}\right.
$$
\end{prop}

{\em Proof}. We shall give only a sketch of the proof.\parskip0pt

By Lemma \ref{hist}(1), (3) we see that
$B_{\alpha+\beta}(y,U,V)\subseteq
B_{\alpha}(x,X,G)\cap B_{\gamma}(x',U,V)$.
Then by Corollary \ref{phar} we conclude
that the set $B_{\alpha+\beta}(y,U,V)$
consists of all elements
$z\in B_{\alpha}(x,X,G)\cap B_{\gamma}(x',U,V)$ such that
the local orbits $V_Uz$ and $V_Uy$ intersect
the same sets from ${\mathcal B}_{<\alpha+\beta}$. \parskip0pt

On the other hand since
$B_{\alpha}(x,X,G)\cap B_{\gamma}(x',U,V)$
is locally $V_U$-invariant, then
$V_{(B_{\alpha}(x,X,G)\cap B_{\gamma}(x',U,V))}z=V_Uz$
whenever $z\in B_{\alpha}(x,X,G)\cap B_{\gamma}(x',U,V)$.
Applying Corollary \ref{phar} to the $G$-space
$\ \langle B_{\alpha}(x,X,G), t^x_{\alpha}\rangle\ $ we conclude
that the set 
$B_{\beta+1}^{t^x_{\alpha}}(y,B_{\alpha}(x,X,G)\cap B_{\gamma}(x',U,V),V)$
consists of all elements $z\in B_{\alpha}(x,X,G)\cap B_{\gamma}(x',U,V)$
such that
the local orbits $V_Uz$ and $V_Uy$ intersect
the same sets from ${\mathcal B}^{t^x_{\alpha}}_{<\beta+1}$.
 
Now the required property follows by induction on $\beta$ with use
of Lemma \ref{hist}(3) both for the original topology and $t^x_{\alpha}$.
$\Box$
\bigskip

This proposition is not involved into main results of the paper. 
For completeness we just describe some application of it.
We start with Hjorth's generalization the notion of a Scott rank.
In \cite{hjorth} Hjorth proves that 
to every $x\in X$ we can assign a cardinal invariant
which can be treated as a counterpat of a Scott rank.
The definition is based on the following lemma.

\begin{lem}\label{height}(Hjorth)
For every $x\in X$ there is some $\gamma <\omega_1$
such that for all $U\in {\mathcal U}$, $V\in \mathcal{V}$ and
$x',x''\in Gx$ we have
$$
(\exists \alpha<\omega_1)\big(B_{\alpha}(x',U,V)\not=B_{\alpha}(x'',U,V )
\big)
\Rightarrow\big(B_{\gamma}(x',U, V)\not=B_{\gamma}(x'', U, V)\big).
$$
\end{lem}

For every $x\in X$ we denote by $\gamma^{\star}(x)$  the least ordinal
$\gamma$ satisfying the statement of Lemma \ref{height}.\parskip0pt

It is proved in \cite{hjorth} that in the case $G=S_{\infty}$
we have $B_{\gamma^{\star}(x)+2}(x,X,G)=Gx$, for every $x\in X$.
It remains true for any closed permutation group but fails
in the general case of an arbitrary Polish group $G$.
Hjorth proves the following weaker statement.

\begin{thm}\label{orb}(Hjorth)
For every $x\in X$,  $U\in {\mathcal U}$, $V\in {\mathcal V}$
and $\alpha\ge \gamma^{\star}(x)+2$ we have
$$
B_{\gamma^{\star}(x)+2}(x,U,V)=B_{\alpha}(x,U,V). 
$$
In particular $B_{\gamma^{\star}(x)+2}(x,X,G)=B_{\alpha}(x,X,G)$.
\end{thm}

The following assertion is a direct consequence of Proposition \ref{subs}.

\begin{cor} For every ordinal $1\le \alpha<\omega_1$ we have
$$
B_{\gamma^{\star}(x)+2}(x,X,G)=
B^{t_x^{\alpha}}_{(\gamma^{t_x^{\alpha}})^{\star}(x)+2}
(x,B_{\alpha}(x,X,G),G).
$$
\end{cor}

\section{Eventually open actions}

The generalized Scott analysis is an important tool
in studying orbit equivalence relations.
The result of Hjorth from \cite{hjorth} which we mentioned in
Introduction can be expressed in terms of $\alpha$-pieces as follows:
\begin{quote}
If the system of generalized Scott invariants for a $G$-space $X$
is complete, i.e. for every $x\in X$ the piece
$B_{\gamma^{\star}(x)+2}(x,X,G)$ coincides with $Gx$,
then the orbit equivalence relation induced by the $G$-action
is classifiable by countable models.
\end{quote}

In this section we study a property of  continuous actions 
of $G$ on $X$ which, as we will see later, is equivalent to the equality
$B_{\gamma^{\star}(x)+2}(x,X,G) =Gx$ for all $x\in X$.

\begin{definicja}
We say that the action $G$ on $X$ is eventually open if
for every $x\in X$ and $V\in {\mathcal V}$ there are $n,m\in\omega$
such that $x\in U_n$ and $(V_m)_{U_n}x\subseteq Vx$.
\end{definicja}

Observe that every Polish group has
an eventually open action, e.g. the action by left multiplication
on itself.
Moreover if a Polish group admits
a basis of open subgroups at its unity then all their
continuous actions are eventually open.
We now show that the converse is not true.
The proof of this assertion is based on the idea described
in \cite{hjorth}.

\begin{prop} All continuous actions of ${\mathbb R}$ on Polish spaces
are eventually open.
\end{prop}

{\em Proof}. Le $X$ be a Polish ${\mathbb R}$-space, $x\in X$ and
$I\subseteq {\mathbb{R}}$ be an open, symmetric interval.
Since $S_x$, the stabilizer of $x$, is a closed subgroup of ${\mathbb R}$,
then either  $S_x={\mathbb{R}}$ or $S_x$ is nowhere dense.
In the first case we have ${\mathbb{R}}x=\{x\}$
and so $J_Ux\subseteq Ix$ for every open $J$ and $U$.
In the second case $(-\infty,0)\setminus S_x\not=\emptyset$
and $(0,\infty)\setminus S_x\not=\emptyset$.
Hence there is $U$ containing $x$ such that both sets
$(-\infty,0)x\setminus \overline{U}$ and $(0,\infty)x\setminus\overline{U}$
are nonempty.
Then by continuity of the action
we can choose $J=(-a,a)\subseteq I$
and $k,l\in {\mathbb{N}}$ such that
$(ka, (k+2)a)x\cap \overline{U}=\emptyset$
and $(-(l+2)a, -la)x\cap \overline{U}=\emptyset$.
Then we see that $J_Ux\subseteq [-la,ka]x$.
We have $\{b\in {\mathbb{R}}:bx\in Ix\}=IS_x$, $IS_x$ is open
and so $[-la,ka]\setminus IS_x$ is a compact subset of ${\mathbb{R}}$.
Since $[-la,ka]x\setminus Ix=([-la,ka]\setminus IS_x)x$,
then $[-la,ka]x\setminus Ix$ is a compact subset of $X$ not containing $x$.
Hence there is a basic open $U_n$ containing $x$ such that
$U_n\subseteq U\setminus ([-la,ka]x\setminus Ix)$.
Since $J_{U_n}x\subseteq J_Ux\cap U_n\subseteq [-la,ka]x\cap U_n$,
then we finally conclude $J_{U_n}x\subseteq Ix$.
$\Box$
\bigskip

Thus for Polish group actions eventual openness is a property weaker
than being induced by a group admitting
a basis of open subgroups at its unity.
Nevertheless as we shall see below eventual openness is equivalent
to completness of the system of generalized Scott invariants.
\parskip0pt

The rest of the section is divided into two parts.
The first one is devoted to local counterparts of Vaught transforms
$\Delta$ and $*$ which will be applied in the second part
in the proof of the result announced above.
Moreover we will see later that eventually open actions are exactly
those for which the map  $G\to Gx$:  $g\to gx$ is open with respect
to $t^x_{\gamma^{\star}(x)+2}$, for every $x\in X$.
This property motivates the name.
At this place note that using Theorem \ref{orb} and
Lemma \ref{hist}(3) we have the following fact.

\begin{prop}\label{open}
Let $G$ be a Polish group, $X$ be a Polish $G$-space and $x\in X$.
The  following statements are equivalent:\parskip0pt

(i)  There is $\alpha\ge 1$ such that
the map $G\to G x$: $g\to g x $
is open with respect to $t^{\alpha}_x$. 

(ii) The map\quad $G\to G x$: $g\to g x $
\quad  is open with respect to $t^x_{\gamma^{\star}(x)+2}$. 
\end{prop} 

\subsection{Local Vaught transforms $\Delta_UV$ and $*_UV$}

In the introductory section we saw that in some special cases local
$V_U$-saturation turns Borel sets into Borel sets.
In general if $A$ is a Borel set, then $V_UA$ is analytic
(recall that $VA$ is analytic whenever $A$ is Borel).
Then the question of some $V_U$-local counterparts of Vaught transforms
arises.

\begin{definicja}\label{transform}
Let $V\in \hat{\mathcal V}$ and $U\subseteq X$ be open.
For every Borel set $A\subseteq X$ we define: 
$$
\begin{array}{l@{\ = \ }l}
A^{{\Delta}_U (V,1)}&(A\cap U)^{\Delta V}\cap U ,\\
A^{{\Delta}_U (V,n+1)}&  (A^{{\Delta}_U (V,n)})^{\Delta V}\cap U ,\\
A^{{\Delta}_U V}&\bigcup\limits_n A^{{\Delta}_U (V,n)} ,
\end{array}
$$
$$
\begin{array}{l@{\ = \ }l}
A^{{*}_U (V,1)}&((A\cap U)\cup (X\setminus U))^{* V}\cap U ,\\
A^{{*}_U (V,n+1)}
&(A^{{*}_U (V,n)} \cup (X\setminus U))^{* V}\cap U ,\\
A^{*_U V}&\bigcap\limits_n A^{*_U(V,n)}.
\end{array}
$$
\end{definicja}

{\bf Remark.} It is clear that if we substitute in the above definition
$A$ by $A\cap U$, then we obtain exactly the same sets, especially
$A^{{\Delta}_U V}=(A\cap U)^{{\Delta}_U V}$ and
$A^{*_U V}=(A\cap U)^{*_U V}$.
Therefore we can limit ourselves to subsets of  $U$,
while discussing properties of $V_U$-local Vaught transforms.\parskip0pt

It is natural to ask if the sequences $(A^{{\Delta}_U (V,n)})$
and $(A^{{*}_U (V,n)})$ are monotone.
The positive answer to this question is one of the consequences
of the following statements.

\begin{lem}\label{nkey}
Let $V\in \hat{\mathcal V}$, $U\subseteq X$ be open,
$x\in U$  and $A\subseteq X$ be a Borel set.
Then for every natural number $n>0$
the following  conditions are satisfied:\parskip0pt

$\begin{array}{l@{\ \mbox{ iff } \ }l}
(1)\ x\in A^{{\Delta}_U (V,n)}&
x\in  A^{\Delta ({\langle V\rangle}^x_U(n))}.\\
(2)\ x\in A^{*_U (V,n)}& x\in A^{* ({\langle V\rangle}^x_U(n))}.
\end{array}$
\end{lem}

{\em Proof}. (1) We proceed by induction.
Since $\{g\in V:gx\in A\cap U\}\subseteq {\langle V\rangle}^x_U(1)$
and ${\langle V\rangle}^x_U(1)$ is open,
then $\{g\in V:gx\in A\cap U\}$ is nonmeager in $V$
if and only if it is nonmeager in ${\langle V\rangle}^x_U(1)$.
This settles the case $n=1$.\parskip0pt

Now assume that  for some $n>0$, every Borel set $B\subseteq X$ and
every $y\in U$ we have
$ y\in B^{{\Delta}_U (V,n)}$ iff
$y\in  B^{\Delta ({\langle V\rangle}^y_U(n))}$.  

If $x\in A^{{\Delta}_U (V,n+1)})$,
then according to Definition \ref{transform},
$x\in U$ and  there is $h\in V$ such that
such that $hx\in A^{{\Delta}_U (V,n)}$.
This implies $h\in {\langle V\rangle}^{x}_U(1)$ and
by the inductive assumption
$hx\in A^{\Delta ({\langle V\rangle}^{hx}_U(n))}$.
Then there is a non-meager set $K\subseteq  {\langle V\rangle}^{hx}_U(n)$
such that $fhx\in A$ whenever $f\in K$.
Hence $Khx\subseteq A$, $Kh$ is a non-meager subset of
${\langle V\rangle}^{hx}_U(n)h$ and
${\langle V\rangle}^{hx}_U(n)h\subseteq_{open}
{\langle V\rangle}^{x}_U(n+1)$.
Therefore $x\in A^{\Delta ({\langle V\rangle}^{x}_U(n+1))}$, which
completes the forward direction.\parskip0pt

To prove the converse suppose that
$x\in A^{{\Delta} ({\langle V\rangle}^x_U(n+1))}$.
This means that the set 
$K=\{f\in G: fx\in A\}$
is nonmeager in ${\langle V\rangle}^x_U(n+1)$.
Then the set
$\hat{K}=\{(g,h): gh\in K\}$
is nonmeager in the set
$\{(g,h): \, h\in {\langle V\rangle}^x_U(1),
\, g\in {\langle V\rangle}^x_U(n),
\, gh\in {\langle V\rangle}^x_U(n+1)\}$.
By Lemma \ref{bH}(3) we see that the latter set
 is an open subset of 
${\langle V\rangle}^x_U(n)\times {\langle V\rangle}^x_U(1)$.
Hence $\hat{K}$ is nonmeager in  the product
${\langle V\rangle}^x_U(n)\times {\langle V\rangle}^x_U(1)$.
On the other hand borelness of $A$ implies borelness of the set
$\hat{K}$.
Thus we can apply Kuratowski-Ulam Theorem to see that
there is a non-meager set $F\subseteq {\langle V\rangle}^x_U(1)$
such that for every $h\in F$ the set $\{g\in G: ghx\in A\}$ is non-meager
in  ${\langle V\rangle}^{hx}_U(n)\}$.
Hence 
$$
x\in (A^{{\Delta}_U (V,n)})^{\Delta V}\cap U
=A^{\Delta ({\langle V\rangle}^x_U(n+1))}. 
$$ 
This completes the proof of (1).
We can prove (2) in a similar way.
$\Box$

\begin{cor} Let $V\in \hat{\mathcal V}$ and $U\subseteq X$ be open.
For every Borel set $A\subseteq X$ the sequence
$(A^{{\Delta}_U (V,n)})$ is increasing while the sequence
$(A^{{*}_U (V,n)})$ is decreasing.
\end{cor}

{\em Proof}. This is an easy consequence of the lemma above,
Definition \ref{H} and standard properties of the original
topological Vaught transforms.$\Box$ 
\bigskip

Let $x\in U$. Intuitively $x$ is an element of $A^{{\Delta}_U V}$
if  $A$ contains some "big" part of its local $V_U$-orbit.
Similarly $x$ belongs to  $A^{*_U V}$ if
"almost whole" $V_Ux$ is contained in $A$.
The following statement is a  precise formulation of this idea.
It also follows directly from Definition \ref{transform} and
Lemma \ref{nkey}.

\begin{cor}\label{key}
Let $V\in \hat{\mathcal V}$, $U\subseteq X$ be open and $A\subseteq X$
be a Borel set.
Then for every $x\in U$ the following conditions are satisfied: 
\parskip0pt

$\begin{array}{l@{\ \mbox{ iff } \ }l}
(1)\ x\in A^{{\Delta}_U V}& x\in A^{\Delta ({\langle V\rangle}_U^x)}.\\
(2)\ x\in A^{*_U V}& x\in A^{* ({\langle V\rangle}_U^x)}.
\end{array}$
\end{cor}

Now it is clear that the basic properties of $V_U$-local Vaught transforms
are similar to the properties of the original (topological) Vaught transforms.

\begin{lem}\label{basic}
Let $V, V'\in \hat{\mathcal V}$, $U\subseteq X$ be open and $x\in U$.
Let $A,A_0,A_1, \ldots , A_n,\ldots $ be Borel subsets of $X$.
Then the following statements hold.
\parskip0pt

$(1)$  If $V'\subseteq V$, then
$A^{{\Delta}_UV'}\subseteq A^{{\Delta}_U V}$ and
$A^{*_UV'}\supseteq A^{*_U V}$.\parskip2pt

$(2)$ $(U\setminus A)^{{\Delta}_U V}=U\setminus A^{{*}_U V}.$\parskip2pt

$(3)$ $(\bigcup\limits_n A_n)^{{\Delta}_U V}=
\bigcup\limits_n (A_n^{{\Delta}_U V})$ and
$(\bigcap\limits_n A_n)^{{*}_U V}=
\bigcap\limits_n (A_n^{{*}_U V}).$\parskip2pt

$(4)$ $x\in A^{{\Delta}_U V}\, \Leftrightarrow
\, x\in \bigcup\{A^{* (V''g)}\cap U:\
g\in  G_0\cap {\langle V\rangle}_U^x,\ V''\subseteq V, V''g\subseteq {\langle V\rangle}_U^x\};$
\parskip2pt

\quad\  $x\in A^{*_U V}\, \Leftrightarrow \,
x\in \bigcap\{A^{\Delta (V''g)}\cap U:
g\in  G_0\cap {\langle V\rangle}_U^x,  V''\subseteq V, V''g\subseteq {\langle V\rangle}_U^x\}$.
\parskip2pt

$(5)$ For every countable ordinal $\alpha>0$ we have:
\parskip0pt

(i) If $A\in \Sigma^0_{\alpha}$,
then $A^{{\Delta}_U V}\in \Sigma^0_{\alpha}$ and
$A^{*_U V}\in \Pi^0_{\alpha+1}$.
\parskip0pt

(ii) If $A\in \Pi^0_{\alpha}$, then $A^{*_U V}\in \Pi^0_{\alpha}$
and  $A^{{\Delta}_U V}\in \Sigma^0_{\alpha+1}$.

\end{lem}

{\em Proof.} (1) and (3) are immediate cosequences
of the analogous properties of the original Vaught transforms. 

(2) By Corollary \ref{key}(1)
$x\in (U\setminus A)^{{\Delta}_U V}$ if and only if
$x\in (U\setminus A)^{\Delta ({\langle V\rangle}_U^x)}$.
By standard properties of the original Vaught transforms, the latter
is equivalent to $x\in U\setminus A^{* {\langle V\rangle}_U^x}$.
Then by Corollary \ref{key}(2) we get
$x\in U\setminus A^{{*}_U V}.$ 

(4) Applying Corollary \ref{key} together
with the properties of the original Vaught transforms we get the equivalence
$$(\star)\ x\in A^{{\Delta}_U V}\ \Leftrightarrow \
x\in \bigcup\{(A^{*W})\cap U: W\subseteq {\langle V\rangle}_U^x
\mbox{ open }\}.$$
Since the group operation is continuous,
for every  basic open  $W\subseteq {\langle V\rangle}_U^x$
we can find $V''\subseteq V$ and $g\in G_0\cap {\langle V\rangle}_U^x$ such that
$V''g\subseteq W$.
Then  we have $A^{*W}\subseteq A^{*(V''g)}$.
This completes the proof of the first equivalence.
We can prove the second one in a similar way. 

(5) We may easily derive it from the definition using induction
and the analogous properties of the original Vaught transforms.
$\Box$
\bigskip

We close the discussion of $V_U$-local Vaught transforms
with the following important property.

\begin{lem}\label{10}  For every open $U\subseteq X$, Borel set
$A\subseteq X$ and
$V\in \hat{{\mathcal V}}$, the following statements hold.
\parskip0pt

$(1)$ $A^{*_U V}\subseteq A^{\Delta _U V}\subseteq V_UA$.\parskip0pt

$(2)$ $A^{{\Delta}_U V}$ and $A^{*_U V}$ are locally
$V_U$-invariant.\parskip0pt

$(3)$ If $A$ is locally $V_U$-invariant, then
$A^{{\Delta}_U V}=A^{*_U V}=A\cap U$.
\end{lem}

{\em Proof.} (1) is immediate by Corollary \ref{key}. 

(2) Accordingly to the remarks following the definition of
local $V_U$-invariantness, we have to prove that
if $x\in A^{\Delta_U V}$ ($x\in A^{*_U V}$)
then $V_Ux\subseteq A^{\Delta_U V}$ (resp. $V_Ux\subseteq A^{*_U V}$). 

Take an arbitrary $x\in A^{\Delta_U V}$ ($x\in A^{*_U V}$).
Then by \ Corollary \ref{key} \ the set 
$K=\{g\in {\langle V\rangle}_U^x:gx\in A\}$
is nonmeager (comeager) in ${\langle V\rangle}_U^x$.
Thus for any $h\in  {\langle V\rangle}_U^x$ the set
$Kh^{-1}$ is nonmeager (comeager) in ${\langle V\rangle}_U^{hx}$.
By Corollary \ref{key} again, the latter  yields
$hx\in A^{{\Delta}_U V}$ (resp. $hx\in A^{*_U V}$).
Thus we are done since  $V_Ux={\langle V\rangle}_U^xx$.
\parskip0pt

(3) is immediate by Corollary \ref{key} and Lemma \ref{10}(1).
$\Box$ 
\bigskip

By point (3) of the lemma above if $A$ is locally $V_U$-invariant,
then $V_UA=A^{{\Delta}_U V}$.
We may generalize the property as follows.

\begin{prop}  Let $V\in \hat{\mathcal V}$, $U\subseteq X$ be open
and  $A\subseteq U$ be a Borel set.
If there are open $U'\subseteq U$ and $V'\in\hat{{\mathcal V}}$ such that
$A\subseteq U'$ and $A$ is locally $V'_{\, U'}$-invariant,
then $V_UA=A^{{\Delta}_U V}$.
\end{prop}

{\em Proof}. We have $\supseteq$ by point (1) of the previous lemma.
\parskip0pt

To prove $\subseteq$ take any element $hx$ where $x\in A$
and $h\in {\langle V\rangle}^x_U$.
Since $A$ is locally $V'_{U'}$-invariant thus
${\langle V'\rangle}^x_Ux\subseteq A$.
Then also ${\langle V'\rangle}^x_Uh^{-1}hx\subseteq A$.
By Lemma \ref{bH} ${\langle V\rangle}_U^{hx}={\langle V\rangle}^x_Uh^{-1}$.
Thus we see that $\{f\in {\langle V\rangle}_U^{hx}:fhx\in A\}$ contains
an open set ${\langle V'\rangle}_{U'}^xh^{-1}$
which proves $hx\in A^{{\Delta}_U V}$.
$\Box$ 
\bigskip

This lemma together with Lemma \ref{basic}(5) yield the following
assertion.

\begin{cor}\label{vi} Let $V\in \hat{\mathcal V}$, $U\subseteq X$ be open,
$A\subseteq U$ and $A\in \Pi^0_{\alpha}(X)$ for some ordinal $\alpha$.
If there are open $U'\subseteq U$ and $V'\in\hat{{\mathcal V}}$ such that
$A\subseteq U'$ and $A$ is locally $V'_{\, U'}$-invariant,
then $V_UA\in \Sigma^0_{\alpha+1}(X)$.
\end{cor}

\subsection{Eventual openness and generalized Scott analysis}

\begin{thm}\label{claschar} Let $G$ be a Polish group and $X$ be
a Polish $G$-space.
The  following statements are equivalent:\parskip0pt

$(1)$ The $G$-action on $X$ is eventually open.\parskip0pt

$(2)$ For every  $V\in {\mathcal V}$, $U\in {\mathcal U}$,
$x\in X$, ordinal $\alpha\ge 1$ and every locally\parskip0pt

$V_U$-invariant Borel set
$A\in \Sigma^0_{\alpha}(X)\cup\Pi^0_{\alpha}(X)$\
if $\ x\in A$, then $B_{\alpha}(x,U,V)\subseteq A.$\parskip0pt

$(3)$ For every $x\in X$ we have
$Gx=B_{\gamma^{\star}(x)+2}(x,X,G)$.\parskip0pt

$(4)$ For every $x\in X$ the map $G\to Gx$:\  $g\to gx$\
is open with respect to $t^x_{\gamma^{\star}(x)+2}$.
\end{thm}

{\em Proof}. (1)$\Rightarrow$(2)
We apply an inductive argument. 

Consider the case $\alpha=1$.
If $A$ is a locally $V_U$-invariant open set containing $x$,
then the local orbit $V_U x$ meets some basic open set $U_n\subseteq A$.
Then $V_U U_n\subseteq A$ and so
$$
A\supseteq \bigcap\{ V_UU_i :U_i\cap V_U x\not=\emptyset\}
\supseteq B_1(x,U,V).
$$ 
If $A$ is a locally
$V_U$-invariant closed set containing $x$, then we have 
$$
A\supseteq\bigcap \{X\setminus V_UU_n: A\cap U_n=\emptyset\}
\supseteq\bigcap \{X\setminus V_UU_n: V_U x\cap U_n=\emptyset\}
\supseteq B_1(x,\sigma). 
$$

To go through the successor step assume that
for every $x'\in Gx$, $U_n$, $V_m$ and every locally
$(V_m)_{U_n}$-invariant $A\in \Sigma^0_{\alpha}\cup \Pi^0_{\alpha}(X)$
we have $B_{\alpha}(x',U_n,V_m)\subseteq A$, whenever $x'\in A$.
Then consider an arbitrary  $V_{U}$-invariant set
$A\in \Sigma^0_{\alpha+1}\cup \Pi^0_{\alpha+1}(X)$
such that $x\in A$.\parskip4pt

$1^o$ If $A\in \Sigma^0_{\alpha+1}(X)$
then  $A$ can be presented as a union $A=\bigcup\limits_i A_i$ such that
$\{A_i:\ i<\omega\}\subseteq \Pi^0_{\alpha}(X)$.
Since $A$ is  locally $V_U$-invariant,
then according to Lemma \ref{10}(3) and Lemma \ref{basic}(3) we have
$A=A^{\Delta_UV}=\bigcup\limits_i A_i^{\Delta_UV}$.
There is $\ i<\omega$ such that $x\in A_i^{\Delta_UV}$.
Then by Lemma \ref{basic}(4),
there are $g\in G_0\cap {\langle V\rangle}_U^x$ and $V'\subseteq V$
such that
$V'g\subseteq {\langle V\rangle}_U^x$ and $x\in A_i^{*(V'g)}$.
This means that
the set $\{h\in G: hx\in A_i\}$ is comeager in $V'g$
or equivalently the set  $T=\{h\in G: hgx\in A_i\}$ is comeager in
$V'$.
Let $G_{gx}$ be the stabilizer of $gx$.
Since $T=TG_{gx}$ and  $V'G_{gx}=\{h\in G: hgx\in V'gx\}$
thus $T$ is comeager in $\{h\in G: hgx\in V'gx\}$.
The $G$-action is eventually open, so there are
$U_n\subseteq U$ containing $gx$ and $V_m\subseteq V'$
such that $(V_m)_{U_n}gx\subseteq V'gx$.
Then $\langle V_m\rangle_{U_n}^{gx}\subseteq \{h\in G: hgx\in V'gx\}$,
and so $T$ is comeager in  ${\langle V_m\rangle}_{U_n}^{gx}$.
By Corollary \ref{key},  the latter implies $gx\in A_i^{*_{U_n}V_m}$.
Using Lemmas \ref{basic}(5) and \ref{10}(3) we see that
$A_i^{*_{U_n}V_m}$ is a locally $(V_m)_{U_n}$-invariant
$\Pi^0_{\alpha}$-set.
Thus we can apply the inductive assumption to get
$B_{\alpha}(gx,U_n,V_m)\subseteq A_i^{*_{U_n}V_m}$.
On the other hand Lemmas  \ref{10}(1) and \ref{basic}(1)
together with the assumption that $A$ is locally $V_U$-invariant imply

\centerline{$A_i^{*_{U_n}V_m}\subseteq A_i^{\Delta_{U_n}V_m}
\subseteq (V_m)_{U_n}A\subseteq A$.}

Therefore $V_UB_{\alpha}(gx,U_n,V_m)\subseteq A$ and so
$B_{\alpha+1}(x,U,V)\subseteq A$.
\parskip4pt

$2^o$ If $A$ is a locally $V_U$-invariant $\Pi^0_{\alpha+1}$-set then
$X\setminus A$ is a locally $V_U$-invariant  $\Sigma^0_{\alpha+1}$-set.
Suppose that $B_{\alpha+1}(x,U,V)\not\subseteq A$.
Then there is $y\in B_{\alpha+1}(x,U,V)$ such that
$y\in X\setminus A$.
By $1^o$ this implies
$B_{\alpha+1}(y,U,V)\subseteq X\setminus A$.
This by Lemma \ref{hist}(1) contradicts the assumption that
$x\in A$.\parskip3pt

(2)$\Rightarrow$(3) Since $Gx$ is an invariant Borel set, then
there is an ordinal $\alpha$ such that
$B_{\alpha}(x,X,G)\subseteq Gx$.
Hence by Theorem \ref{orb} we have
$B_{\gamma^{\star}(x)+2}(x,X,G)\subseteq Gx$.
\parskip3pt

(3)$\Rightarrow$(4) By Effros theorem on $G_{\delta}$-orbits.
\parskip3pt

(4)$\Rightarrow$(1) Take an arbitrary $V\in {\mathcal V}$.
The set $Vx$\ is  \ $t^x_{\gamma^{\star}(x)+2}$-open,
so there are $n,m\in \omega$ such that
$B_{\gamma^{\star}(x)+1}(x,U_n,V_m)\cap Gx
\subseteq Vx$.
Hence $(V_m)_{\, U_n}x\subseteq Vx$.
$\Box$
\bigskip

If $X$ is a Polish $G$-space under an eventually open action,
then by the theorem of Hjorth (see Lemma 6.30. \cite{hjorth})
the orbit equivalence relations iduced on $X$
is classifiable by countable models.
Thus for every ordinal $\alpha\ge 1$ the orbit equivalence relation
induced on the Polish $G$-space
$\langle B_{\alpha}(x,X,G), t^x_{\alpha}\rangle$ is also classifiable
by countable models. By another theorem of Hjorth
(see Corollary 3.19. \cite{hjorth}) this in turn implies that the $G$-action
on $B_{\alpha}(x,X,G)$  is not turbulent
with respect to $t^x_{\alpha}$.
In particular the action of $G$ on $B_{\gamma^{\star}(x)+2}(x,X,G)$
is not turbulent with respect to
$t^x_{\gamma^{\star}(x)+2}$.

%
%
%
\end{document}